\documentclass[10pt]{amsart}
\calclayout
\usepackage[british]{babel}
\usepackage[utf8]{inputenc}
\usepackage{amsfonts}
\usepackage{amssymb}
\usepackage{amsmath}
\usepackage{amsthm}
\usepackage{bm}
\usepackage{latexsym}
\usepackage{amscd}
\usepackage{mathrsfs}
\usepackage{graphicx}
\usepackage{dsfont}
\usepackage{enumitem}
\usepackage{mathtools}
\usepackage{enumerate}
\usepackage{url}
\usepackage[dvipsnames]{xcolor}
\usepackage{hyperref}
\hypersetup{colorlinks=true,linkcolor=ForestGreen,citecolor=Orange,urlcolor=blue}
\usepackage[capitalise]{cleveref}
\usepackage{tikz}
\usetikzlibrary{patterns,patterns.meta}
\tikzdeclarepattern{
  name=mydlines,
  parameters={
    \pgfkeysvalueof{/pgf/pattern keys/size},
    \pgfkeysvalueof{/pgf/pattern keys/angle},
    \pgfkeysvalueof{/pgf/pattern keys/line width},
  },
  bounding box={
    (0,-0.5*\pgfkeysvalueof{/pgf/pattern keys/line width}) and
    (\pgfkeysvalueof{/pgf/pattern keys/size},
    0.5*\pgfkeysvalueof{/pgf/pattern keys/line width})},
  tile size={(\pgfkeysvalueof{/pgf/pattern keys/size},
    \pgfkeysvalueof{/pgf/pattern keys/size})},
  tile transformation={rotate=\pgfkeysvalueof{/pgf/pattern keys/angle}},
  defaults={
    size/.initial=5pt,
    angle/.initial=45,
    line width/.initial=.4pt,
  },
  code={
    \draw [dashed,
           line width=\pgfkeysvalueof{/pgf/pattern keys/line width}]
    (0,0) -- (\pgfkeysvalueof{/pgf/pattern keys/size},0);
  },
}
\tikzdeclarepattern{
  name=myslines,
  parameters={
    \pgfkeysvalueof{/pgf/pattern keys/size},
    \pgfkeysvalueof{/pgf/pattern keys/angle},
    \pgfkeysvalueof{/pgf/pattern keys/line width},
  },
  bounding box={
    (0,-0.5*\pgfkeysvalueof{/pgf/pattern keys/line width}) and
    (\pgfkeysvalueof{/pgf/pattern keys/size},
    0.5*\pgfkeysvalueof{/pgf/pattern keys/line width})},
  tile size={(\pgfkeysvalueof{/pgf/pattern keys/size},
    \pgfkeysvalueof{/pgf/pattern keys/size})},
  tile transformation={rotate=\pgfkeysvalueof{/pgf/pattern keys/angle}},
  defaults={
    size/.initial=5pt,
    angle/.initial=45,
    line width/.initial=.4pt,
  },
  code={
    \draw [solid,
           line width=\pgfkeysvalueof{/pgf/pattern keys/line width}]
    (0,0) -- (\pgfkeysvalueof{/pgf/pattern keys/size},0);
  },
}
\usepackage{pgfplots}
 \pgfplotsset{
     tick align=outside,
     x grid style={white},
     xmajorgrids,
     y grid style={white},
     ymajorgrids,
     axis line style={white},
     axis background/.style={fill=white!92!black},
     legend style={draw=white, fill=white},
     legend cell align={left}
 }

\numberwithin{equation}{section}
\newtheorem{theorem}{Theorem}
\newtheorem{lemma}[theorem]{Lemma}

\newtheorem{remark}[theorem]{Remark}

\newtheorem{hypothesis}{Hypothesis}

\theoremstyle{definition}

\newcommand{\ee}{\mathrm{e}} 

\newcommand{\dd}{\, \mathrm{d}} 
\newcommand{\var}{\varepsilon}

\newcommand{\N}{{\mathbb N}}

\newcommand{\R}{{\mathbb R}}
\newcommand{\T}{{\mathbb T}}

\renewcommand{\S}{{\mathbb S}}

\newcommand{\sD}{\mathsf D}

\newcommand{\sU}{\mathsf U}

\newcommand{\cC}{\mathcal C}
\newcommand{\cD}{\mathcal D}

\newcommand{\weakto}{\rightharpoonup}

\newcommand{\vAvg}[1]{\langle{#1}\rangle} 
\newcommand{\gAvg}[1]{{\langle\!\langle}{#1}{\rangle\!\rangle}} 

\newcommand{\supp}{\operatorname{supp}} 
\newcommand{\indicator}{\mathbf{1}} 
\newcommand{\init}{\mathrm{in}} 
\newcommand{\LB}{\mathrm{LB}} 

\newcommand{\hypref}[1]{\hyperref[#1]{\textbf{(H\ref*{#1})}}}

\renewcommand{\vec}[1]{\mathbf{#1}}

\begin{document}

\title[Trajectorial hypocoercivity]{Trajectorial
  hypocoercivity and application to control theory}

\author[H. Dietert]{Helge Dietert} \address{Université Paris Cité and
  Sorbonne Université, CNRS\\ IMJ-PRG, F-75006 Paris, France.}
\email{helge.dietert@imj-prg.fr}

\author[Hérau]{Frédéric Hérau} \address{Laboratoire de Mathématiques
  Jean Leray, Nantes Université\\ 2 rue de la Houssinière, BP 92208
  F-44322 Nantes Cedex 3, France}
\email{frederic.herau@univ-nantes.fr}

\author[H. Hutridurga]{Harsha Hutridurga} \address{Department of
  Mathematics, Indian Institute of Technology Bombay, Powai,
  Mumbai 400076, India.}  \email{hutri@math.iitb.ac.in}

\author[C. Mouhot]{Clément Mouhot} \address{Department of Pure
  Mathematics and Mathematical Statistics, University of
  Cambridge, Wilberforce Road, CB3 0WA Cambridge, UK}
\email{c.mouhot@dpmms.cam.ac.uk}

\date{\today}

\begin{abstract}
  We present the quantitative method of the recent  work~\cite{dietert-herau-hutridurga-mouhot-2022-preprint-quantitative-geometric-control-linear-kinetic-theory} in a simple setting, together with a compactness argument that was not included in~\cite{dietert-herau-hutridurga-mouhot-2022-preprint-quantitative-geometric-control-linear-kinetic-theory} and has interest per se. We are concerned with the exponential stabilisation (spectral gap) for linear kinetic equations with degenerate thermalisation, i.e. when the collision operator vanishes on parts of the spatial domain. The method in~\cite{dietert-herau-hutridurga-mouhot-2022-preprint-quantitative-geometric-control-linear-kinetic-theory} covers both scattering and Fokker-Planck type operators, and deals with external potential and boundary conditions, but in these notes we present only its core argument and restrict ourselves to the kinetic Fokker-Planck in the periodic torus with unit velocities and a thermalisation degeneracy (this equation is not covered by the previous results~\cite{bernard-salvarani-2013-degenerate-linear-boltzmann,han-kwan-leautaud-2015-geometric,evans-moyano-2019-preprint-quantitative-rates-convergence-equilibrium-degenerate}).
\end{abstract}



\maketitle
\tableofcontents

\section{Introduction}

\subsection{The setting}

Let us consider the linear kinetic equation
\begin{equation}
  \label{eq:evo}
  \partial_t f + v \cdot \nabla_x f = \sigma \Delta_{\LB} f
\end{equation}
for a time-dependent probability density \(f(t,x,v) = f_t(x,v)\) over
the phase space \((x,v) \in \T^d \times \S^{d-1}\), where \(\T^d\) is
the unit torus and \(\S^{d-1} = \{v \in \R^d : |v| = 1\}\), modelling
massless or nearly massless particles with unit velocities. The right
hand side \(\Delta_{\LB}\) is the Laplace-Beltrami operator and
corresponds to the classical Fokker-Planck operator when velocities
are restricted to the sphere. Finally, \(\sigma = \sigma(x) \in L^\infty(\T^d;[0,\infty))\) is a weight that can vanish in part of the spatial domain and models the thermalisation degeneracy.

The evolution~\eqref{eq:evo} has the stationary state \(f_\infty = 1\) and conserves the mass \(\int f \dd x \dd v\). When \(\sigma \equiv 1\), it is one of the simplest examples of hypocoercive equation: any solution \(f_t\) associated to the initial data \(f_\init \in L^2 \) with zero mass \(\int f_\init \dd x \dd v=0\) will converge exponentially to zero. A natural question arises then, inspired from control theory: under what conditions on \(\sigma\), will the evolution~\eqref{eq:evo} yield exponential relaxation to equilibrium?

When the right hand side in~\eqref{eq:evo} is a bounded integral
scattering operator (linear Boltzmann or relaxation operator), this
question has been answered by
\cite{bernard-salvarani-2013-degenerate-linear-boltzmann,han-kwan-leautaud-2015-geometric}
by compactness arguments when $\sigma$ satisfies a \emph{geometric
  control condition} borrowed from control theory of wave
equations~\cite{MR1178650}.  However these works crucially rely on the
facts that (1) the right hand side operator is bounded, and (2) writes
as a non-negative integral operator minus a local part. The case we
consider here is conceptually different, and requires new
methods. Another direction for bounded operators is given in
\cite{evans-moyano-2019-preprint-quantitative-rates-convergence-equilibrium-degenerate}
who answered quantitatively by Harris theorem from probability theory.

\subsection{The geometric control condition}
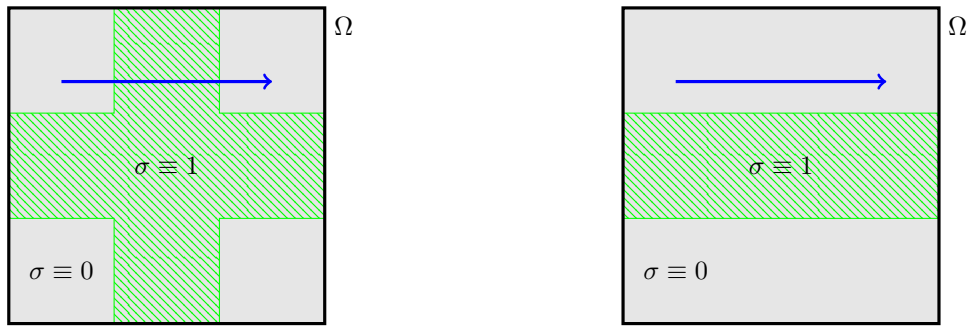
\begin{figure}[htb]
  \begin{centering}
    \begin{tikzpicture}[scale=1.4]
      \draw[fill, color=black!10!white] (0,0) -- (3,0) -- (3,3) -- (0,3) -- cycle;
      \draw[green, pattern=north west lines, pattern color=green]
      (0,1) -- (0,2) -- (1,2) -- (1,3) -- (2,3) -- (2,2) -- (3,2) --
      (3,1) -- (2,1) -- (2,0) -- (1,0) -- (1,1) -- cycle;
      \draw[very thick] (0,0) -- (3,0) -- (3,3) -- (0,3) -- cycle;
      \draw (3,3) node[anchor=north west] {$\Omega$};
      \draw (1.5,1.5) node {$\sigma \equiv 1$};
      \draw (0.5,0.5) node {$\sigma \equiv 0$};
      \draw[blue,very thick,->] (0.5,2.3) -> (2.5,2.3);
    \end{tikzpicture}
    \hfill
    \begin{tikzpicture}[scale=1.4]
      \draw[fill, color=black!10!white] (0,0) -- (3,0) -- (3,3) -- (0,3) -- cycle;
      \draw[green, pattern=north west lines, pattern color=green]
      (0,1) -- (0,2) -- (3,2) -- (3,1) -- cycle;
      \draw[very thick] (0,0) -- (3,0) -- (3,3) -- (0,3) -- cycle;
      \draw (3,3) node[anchor=north west] {$\Omega$};
      \draw (1.5,1.5) node {$\sigma \equiv 1$};
      \draw (0.5,0.5) node {$\sigma \equiv 0$};
      \draw[blue,very thick,->] (0.5,2.3) -> (2.5,2.3);
    \end{tikzpicture}
  \end{centering}
  \caption{Illustration of \hypref{h:gcc}. On the left, all lines hit
    and spend a controlled fraction of time in the thermalisation set
    $\Sigma = \supp \sigma$ on a time interval $[0,T]$, yielding
    GCC and exponential convergence. On the right, there is a set of
    configurations with zero measure whose trajectories never hit
    $\Sigma$, and around this set the time to hit $\Sigma$ can be
    arbitrarily large: only the non-uniform GCC holds, and one
    typically expects polynomial rate of convergence.}
  \label{fig:gcc}
\end{figure}

The transport equation on the left hand side of~\eqref{eq:evo} is
solved by the characteristics
\begin{equation*}
  Z_t(x,v):=(X_t(x,v),V_t(x,v)):=(x+tv,v).
\end{equation*}
The (uniform) \emph{geometric control condition} (GCC) intuitively
means that, in a given fixed time, all trajectories spends a positive
time (bounded below) in a region where $\sigma \gtrsim 1$ (where
thermalisation truly occurs), see Figure~\ref{fig:gcc}. The
non-uniform version of this condition intuitively means that all
trajectories eventually enter the support of $\sigma$ (without
restricting the time horizon or asking that the trajectories spend
time in a region where $\sigma$ remains strictly away from zero).

We adopt the following precise definition:
\begin{hypothesis}[Geometric control condition]
  \label{h:gcc}
  The (uniform) GCC writes
  \begin{equation}
    \label{eq:gcc-cont}
    \exists \, T^*,c>0 \ \text{ such that } \ \forall \, (x,v) \in \T^d \times \S^{d-1}:  \quad \int_0 ^{T^*} \sigma(X_t(x,v)) \dd t \ge c.
  \end{equation}
  We also assume $\sigma \in C^1$, and therefore~\eqref{eq:gcc-cont} implies that there is \(\Sigma \subset \T^d\) open with \(\cC^1\) boundary and a smooth \(\chi : \T^d \to [0,\infty)\), so that \(\indicator_{\Sigma} \lesssim \sigma\), \(\supp \chi \subset \Sigma\) and
  \begin{equation}
    \label{eq:gcc}
    \forall \, (x,v) \in \T^d \times \S^{d-1}:
    \quad
    \int_0 ^{T^*} \chi(X_t(x,v)) \dd t \ge 1.
  \end{equation}
\end{hypothesis}

\begin{remark}
  In this control condition, we assume slightly more regularity on
  $\sigma$ than in the literature:
  \cite{bernard-salvarani-2013-degenerate-linear-boltzmann} only needs
  \(\sigma \in L^\infty\) while
  \cite{han-kwan-leautaud-2015-geometric} assumes that \(\sigma\)
  continuous. The regularity is only used to conclude \eqref{eq:gcc}
  and we did not try to optimise this assumption.
\end{remark}

Our main result is:
\begin{theorem}[Exponential stabilization]
  Assume \(\sigma\) is bounded and satisfies \hypref{h:gcc} with $\Sigma$ having finitely many connected components. Then there are \(C \ge 1\) and \(\Lambda>0\) such that for any initial data \(f_\init \in L^2\) the corresponding solution \(f_t\) to~\eqref{eq:evo} satisfies
  \begin{equation*}
    \left\| f_t - \left( \int_{\T^d \times \S^{d-1}} f_\init \right) \right\|_{L^2(\T^d \times \S^{d-1})} \le C \ee^{-\Lambda t} \left\| f_\init - \left( \int_{\T^d \times \S^{d-1}} f_\init \right) \right\|_{L^2(\T^d \times \S^{d-1})}.
  \end{equation*}
  Such $C,\Lambda$ can be computed from the proof and only
  depend on $T^*,c,\Sigma,\|\chi\|_{W^{1,\infty}},\|\sigma\|_\infty$.
\end{theorem}

\begin{remark}
  The  \emph{non-uniform GCC} means that for almost every $(x,v) \in \T^d \times \S^{d-1}$
  \begin{equation}
    \label{eq:ngcc-cont}
    \exists \, T=T(x,v) >0 \ \text{ such that } \ \int_0 ^{T(x,v)} \sigma(X_t(x,v)) \dd t >0.
  \end{equation}
  When one replaces~\eqref{eq:gcc-cont} by the non-uniform
  condition~\eqref{eq:ngcc-cont}, our method can be used to prove the
  convergence
  $f_t \to f_\infty = ( \int_{\T^d \times \S^{d-1}} f_\init)$,
  although without a rate.
\end{remark}

\section{Trajectorial approach to hypocoercivity}

\subsection{Fixing the global average}

Consider $f_{\init} \in L^2$ and its associated solution $f_t$. Since the mass is conserved and the equation is linear, $g_t:= f_t - (\int f_\init(x,v) \dd x \dd v)$ is solution to~\eqref{eq:evo} with zero mass
\begin{equation}
  \label{eq:zeromass}
  \int_{\T^d \times \S^{d-1}} g_\init(z) \dd z =0,
\end{equation}
and therefore its associated equilibrium is zero.

\subsection{The local projection}

The local equilibrium is $M(v) = |\S^{d-1}|^{-1}$ and we define
the spatial density (velocity average)
\begin{equation}
  \label{eq:def:density}
  \vAvg{g}(t,x) := \int_{\S^{d-1}} g(t,x,v) \dd v.
\end{equation}

\subsection{The energy estimate}

The $L^2$ norm is the natural entropy for this linear model, and the $H$ theorem takes the form of the energy estimate
\begin{equation}
  \label{eq:dissipation}
  \cD(g_t) := -\frac{\mathrm{d}}{\mathrm{d}t} \| g_t \|_{L^2(\T^{d} \times \S^{d-1})}^2 = \int_{\T^d \times \S^{d-1}} \sigma |\nabla_v g_t|^2 \dd x \dd v \ge 0
\end{equation}
where the gradient in $v$ is the differential tangential to the unit sphere.

\subsection{Integral criterion for exponential stabilization}

We first prove a simple sufficient time-integrated entropy production
inequality that implies exponential convergence. Such criterion is
standard in kinetic theory, and (at least) already appears in a
compactness argument in \cite{Guo-2002-I}.  Exponential decay holds if
and only if there are $T>0$, $\lambda>1$ such that
\begin{equation}
  \label{eq:sufficient}
  \| g_{\init} \|_{L^2(\T^d\times \S^{d-1})}^2 \le \lambda \int_0^{T} \cD(g_t) \dd t.
\end{equation}
More precisely: when~\eqref{eq:sufficient} holds, then
\begin{equation*}
  \| g_t \|_{L^2(\T^d \times \S^{d-1})}
  \le
  \sqrt{\frac{\lambda}{\lambda-1}}
  \exp\left(
    - \frac 1T \log\left(\frac{\lambda}{\lambda-1}\right)
    t
  \right)
  \| g_\init \|_{L^2(\T^d \times \S^{d-1})}.
\end{equation*}

\subsection{Micro-coercivity}

The Poincaré inequality holds in the compact smooth manifold
$\S^{d-1}$: there is $C_P>0$ so that for any $x \in \T^d$
\begin{equation}
  \label{eq:micro}
   \int_{\S^{d-1}} \left|g_t - \vAvg{g_t} M \right|^2 \dd v \le C_P \int_{\S^{d-1}} \sigma |\nabla_v g_t|^2 \dd v.
\end{equation}
This provides control over $(g-\vAvg{g}M)$ on
$\supp \sigma \times \S^{d-1}$ hence on the good set
$\Sigma \times \S^{d-1}$.

\subsection{Following the characteristics to transfer the control}

The next step is to transfer the control of $(g-\vAvg{g} M )$ on the
good set to the whole domain, by following trajectories. Let us prove
that there are $C_1,C_2>0$ so that
\begin{equation}
  \label{eq:following}
  \| g_{\init} \|_{L^2(\T^d\times \S^{d-1})}^2 \le C_1 \int_0^{T^*} \cD(g_t) \dd t + C_2 \int_0^{T^*} \int_{\Sigma} \vAvg{g_t}^2 \dd x \dd t.
\end{equation}
To prove this first write the evolution equation for $g^2$:
\begin{equation*}
  \partial_t \left( g^2_t \right) + v \cdot \nabla_x \left( g^2_t \right) = 2 \left( \Delta_{\LB} g_t \right) g_t
\end{equation*}
and second write it in Duhamel form along the transport flow (writing
$z:=(x,v)$)
\begin{equation*}
  g_t(z)^2 = g_\init(Z_{-t}(z))^2
  + 2 \int_0^t \sigma(Z_{t-s}(z))\, (\Delta_{\LB} g_s)(Z_{t-s}(z))\, g_s(Z_{t-s}(z)) \dd s
\end{equation*}
and third integrate it against $\chi$ from~\hypref{h:gcc} on $[0,T^*] \times \T^d \times \S^{d-1}$:
\begin{multline*}
  \int_0 ^{T^*} \int_{\T^d \times \S^{d-1}} g_t^2(z) \chi(z) \dd t \dd z = \int_{\T^d \times \S^{d-1}} g_{\init}(z)^2 \left( \int_0 ^{T^*} \chi(Z_t(z)) \dd t \right) \dd z \\
  + 2 \int_0 ^{T^*} \int_0^t \int_{\T^d \times \S^{d-1}} \sigma(z)\,
  ( \Delta_{\LB} g_s)(z)\, g_s(z) \chi(Z_{t-s}(z)) \dd z \dd s \dd t
\end{multline*}
where we have used the unitary change of variables $z \mapsto Z_t(z)$ and $z \mapsto Z_{t-s}(z)$.

Now observe that~\eqref{eq:gcc} in~\hypref{h:gcc} implies
\begin{equation*}
  \int_{\T^d \times \S^{d-1}} g_{\init}(z)^2 \left( \int_0 ^{T^*} \chi(Z_t(z)) \dd t \right) \dd z \ge \| g_\init \|^2_{L^2(\T^d \times \S^{d-1})}
\end{equation*}
and $\supp \sigma \subset \Sigma$ implies
\begin{equation*}
  \int_0 ^{T^*} \int_{\T^d \times \S^{d-1}} g_t^2(z) \chi(z) \dd t \dd z \le \| \chi \|_\infty \int_0 ^{T^*} \| g_t \|^2_{L^2(\Sigma \times \S^{d-1})} \dd t.
\end{equation*}
As for the last term we perform an integration by parts:
\begin{multline*}
  2 \int_0 ^{T^*} \int_0^t \int_{\T^d \times \S^{d-1}} \sigma(z)
  \left( \Delta_{\LB} g_s \right)(z) g_s(z) \chi(Z_{t-s}(z)) \dd z \dd
  s \dd t = \\
  - 2 \int_0 ^{T^*} \int_0^t \int_{\T^d \times \S^{d-1}} \sigma(z)
  \left| \nabla_v g_s(z)\right|^2 \chi(Z_{t-s}(z)) \dd z \dd s \dd t
  \\
  - 2 \int_0 ^{T^*} \int_0^t (t-s) \int_{\T^d \times \S^{d-1}}
  \sigma(z) \nabla_v g_s(z) \cdot (\nabla \chi)(Z_{t-s}(z)) g_s(z) \dd
  z \dd s \dd t \\
  \ge - T^*
  \left( 2 \| \chi \|_\infty
    + \frac{1}{\var} \| \nabla \chi \|_\infty^2 \| \sigma\|_\infty
  \right)
  \int_0 ^{T^*} \cD(g_s) \dd s
  - \var T^* \int_0 ^{T^*} \| g_s \|^2 _{L^2(\T^d \times
    \S^{d-1})} \dd s \\
  \ge - T^*
  \left( 2 \| \chi \|_\infty
    + \frac{1}{\var} \| \nabla \chi \|_\infty^2 \| \sigma\|_\infty
  \right)
  \int_0 ^{T^*} \cD(g_s) \dd s
  - \var \left(T^*\right)^2 \| g_\init \|^2 _{L^2(\T^d \times \S^{d-1})}
\end{multline*}
where in the last line we have used that the $L^2$ norm is
non-increasing. Therefore we deduce
\begin{multline*}
  \| g_\init \|^2_{L^2(\T^d \times \S^{d-1})} \le \| \chi \|_\infty \int_0 ^{T^*} \| g_t \|^2_{L^2(\Sigma \times \S^{d-1})} \dd t \\ + T^* \left( 2 \| \chi \|_\infty + \frac{1}{\var} \| \nabla \chi \|_\infty \right) \int_0 ^{T^*} \cD(g_s) \dd s +\var \left(T^*\right)^2 \| \nabla \chi \|_\infty \| g_\init \|^2 _{L^2(\T^d \times \S^{d-1})}.
\end{multline*}
We now use the micro-coercivity~\eqref{eq:micro}:
\begin{align*}
  \int_0 ^{T^*} \| g_t \|^2_{L^2(\Sigma \times \S^{d-1})} \dd t
  &\le \frac{2}{|\S^{d-1}|} \int_0 ^{T^*} \| \vAvg{g_t} \|^2_{L^2(\Sigma)} \dd t + 2 \int_0 ^{T^*} \| g_t - \vAvg{g_t} M \|^2_{L^2(\Sigma \times \S^{d-1})} \dd t \\
  &\le \frac{2}{|\S^{d-1}|} \int_0 ^{T^*} \| \vAvg{g_t} \|^2_{L^2(\Sigma)} \dd t + 2 C_P \int_0 ^{T^*} \cD(g_t) \dd t.
\end{align*}
Taking $\var = (T^*)^{-2} \| \nabla \chi \|_\infty^{-1}/2$, we finally deduce
\begin{multline*}
  \| g_\init \|^2_{L^2(\T^d \times \S^{d-1})} \le \frac{4 \| \chi \|_\infty}{|\S^{d-1}|} \int_0 ^{T^*} \| \vAvg{g_t} \|^2_{L^2(\Sigma)} \dd t \\
  + \left[ 4 C_P \| \chi \|_\infty + 4T^* \| \chi\|_\infty
    + 4 \left( T^* \right)^3 \| \nabla \chi \|_\infty ^2
    \|\sigma\|_\infty \right] \int_0 ^{T^*} \cD(g_t) \dd t
\end{multline*}
which proves~\eqref{eq:following}. We are left with the control of the
local projection on the good set, which is the object of the next two
sections.

\section{The compactness argument}

Assume~\eqref{eq:sufficient} to be false with $T=T^*$: there is  a contradiction sequence \((g^n)_{n\in\N}\) of
solutions with initial data \(g^n_\init\) such that \(\|g^n_\init\| = 1\) (normalised by linearity) and
\begin{equation*}
  \int_0^{T^*} \cD(g^n_t) \dd t \to 0.
\end{equation*}
By weak compactness we then find a subsequence such that
\begin{equation*}
  g^{n'} \weakto g^* \text{ in } L^2([0,T^*] \times \T^d \times \S^{d-1})
  \quad\text{and}\quad
  g^{n'}_\init \weakto g^*_\init \text{ in } L^2(\T^d \times \S^{d-1}).
\end{equation*}
Moreover the velocity averaging lemma ensures that $\vAvg{g^n}$ is
relatively compact for the strong topology in
$L^2([0,T^*] \times \T^d)$. Therefore, we can furthermore assume that our
subsequence satisfies $\vAvg{g^{n'}} \to \vAvg{g^*}$ strongly in
$L^2([0,T^*] \times \T^d)$.

The limit then satisfies \(\partial_t g^* + v \cdot \nabla_x g^* = 0\)
in the weak sense and \(\int_0^{T^*}\cD(g_t^*) \dd t = 0\). This
implies that $g^*$ is constant in $\Sigma \times \S^{d-1}$ since it
has to be constant along the transport flow and equal to its velocity
average.  By connecting any point to a point in
$\Sigma \times \S^{d-1}$ (using the GCC), we deduce that $g^*$ is
constant in $\T^d \times \S^{d-1}$. The weak \(L^2\) convergence implies that
$\int_{[0,T^*] \times T^d \times \S^{d-1}} g^* = \lim_{n'\to \infty}
\int_{[0,T^*] \times T^d \times \S^{d-1}} g^{n'} = 0$ (recall that we
have set the total mass of each $g_n$ to zero). Therefore
$g^* \equiv 0$.

Using the strong convergence of the velocity average
$\vAvg{g} \to \vAvg{g^*}$ in
\(L^2([0,T^*] \times \T^d \times \S^{d-1})\) and taking the limit in
\eqref{eq:following}, we deduce that
\(\| \vAvg{g^*} \|_{L^2([0,T^* \times \Sigma \times \S^{d-1})} \gtrsim
1\). This contradicts $g^* \equiv 0$ proved in the previous
paragraph, and \eqref{eq:sufficient} is proved.

\section{Getting quantitative: the divergence inequality}

We now replace the previous non-constructive argument based on compactness and contradiction by a quantitative one. For the sake of readability, we first assume that the \(\Sigma\) from \hypref{h:gcc} is connected and explain how this can be relaxed at the end.

In view of~\eqref{eq:sufficient} and~\eqref{eq:following} and the fact
that the $L^2$ norm is non-increasing, to close a complete
quantitative argument it is enough to prove that for any $\delta>0$
there is $C_\delta>0$ so that
\begin{equation}
  \label{eq:quant}
  \int_0^{T^*} \int_{\Sigma} \vAvg{g_t}^2 \dd x \dd t \le C_\delta \int_0^{T^*} \cD(g_t) \dd t + \delta \int_0 ^{T^*} \| g_t \|^2_{L^2(\T^d \times \S^{d-1})} \dd t.
\end{equation}

We then define the global average over the good set as
\begin{align*}
  \gAvg{g} := \frac{1}{m} \int_{[0,T^*] \times \Sigma} \vAvg{g}(t,x)
             \dd t \dd x
  \quad \mbox{ with } m & := |[0,T^*] \times \Sigma|
\end{align*}
and then split the term to estimate as
\begin{align*}
  \int_0^{T^*} \int_{\Sigma} \vAvg{g_t}^2 \dd x \dd t
  \le \int_0^{T^*} \int_{\Sigma} \left( \vAvg{g_t} - \gAvg{g}
  \right)^2 \dd x \dd t
  + m \gAvg{g} ^2 =: I_1 + I_2.
\end{align*}

To control $I_1$ we use the following result that goes back
to~\cite{MR553920,MR631691,bourgain-brezis-2002-divergence-equation}:
\begin{lemma}[Divergence inequality]
  \label{lem:div}
  Given $\sU \subset \R^n$, $n \ge 1$, an open connected bounded $C^1$ domain, there is $C_\sD > 0$ and a linear map $\sD$ mapping any $h \in L^2(\sU)$ with $\int_{\sU} h =0$ to a $\vec{F} : \sU \to \R^n$ in $H^1(\sU)$ that satisfies
  \begin{align}
    \label{eq:divergence}
    \left\{
    \begin{aligned}
      &\nabla \cdot \vec{F} = h \text{ in } \sU, \\[2mm]
      &\vec{F} = 0 \text{ on } \partial \sU, \\[2mm]
      &\| \vec{F} \|_{H^1(\sU)} \le C_{\sD} \| h \|_{L^2(\sU)}.
      \end{aligned}
        \right.
  \end{align}
\end{lemma}
This is proved constructively
in~\cite{MR553920,MR631691,bourgain-brezis-2002-divergence-equation},
and we refer to our full
paper~\cite{dietert-herau-hutridurga-mouhot-2022-preprint-quantitative-geometric-control-linear-kinetic-theory}
for extensions of this result to general domains with external
potentials and boundary conditions.

We apply \cref{lem:div} to $h:= \vAvg{g_t} - \gAvg{g}$ on
$\sU := (0,T^*) \times \Sigma$ (with zero mass): there is
$\vec{F} \in H^1(\sU)$ so that~\eqref{eq:divergence} holds, and we
write (using the Dirichlet conditions)
\begin{align}
  \nonumber
  \int_0^{T^*} \int_{\Sigma} \left( \vAvg{g_t} - \gAvg{g} \right)^2 \dd x \dd t
  & = \int_\sU \left( \vAvg{g_t} - \gAvg{g} \right) \left( \nabla_{t,x} \cdot \vec{F} \right)  \dd x \dd t \\
  \label{eq:inverse-div}
  & = - \int_\sU \vec{F} \cdot \nabla_{t,x} \vAvg{g_t} \dd x \dd t.
\end{align}

Denote $\partial_0 = \partial_t$ and $\partial_i = \partial_{x_i}$ for
$i =1,\dots,d$. Then there is $C_3 >0$ so that
\begin{align}
  \label{eq:claim-dec}
   & \forall \, i=0,\dots,d, \quad
    \partial_i \langle g_t \rangle = K_i + \sum_{j=0} ^d
  \partial_j J_{ij} \quad \text{ with}\\
  \label{eq:claim-am}
   & \sum_{i=0} ^d \| K_i \|_{L^2(\sU)} ^2 + \sum_{i,j=0} ^d \|J_{ij}\|_{L^2(\sU)} ^2 \le C_3 \int_0^{T^*} \cD(g_t) \dd t.
\end{align}
Indeed define
$\varphi_i \in C^2(\S^{d-1})$, $i=0,\dots,d$ so that, denoting
$v_0=1$,
\begin{align*}
  & \int_{\S^{d-1}} \varphi_i(v) v_j \dd v = \delta_{ij}, \ i,j=0,\dots,d, \quad \text{ and so} \\
  & \int_{\S^{d-1}} \Big\{ \left( \partial_t + v \cdot \nabla_x \right) \big[ \langle g_t \rangle \big] \Big\} \varphi_i \dd v = \partial_i \langle g_t \rangle.
\end{align*}
The evolution equation on $g$ then implies
\begin{align*}
  \partial_i \langle g_t \rangle
  & = \sigma \int_{\S^{d-1}} \left( \Delta_\LB g_t \right) \varphi_i \dd v + \int_{\S^{d-1}} \Big\{ \left( \partial_t + v \cdot \nabla_x \right) \big[ \langle g_t \rangle M - g_t \big] \Big\} \varphi_i \dd v \\
  & = K_i + \sum_{j=0} ^d \partial_j J_{ij}
\end{align*}
with
\begin{align*}
  \begin{dcases}
    K_i(t,x) := \sigma \int_{\S^{d-1}} \left( g_t - \vAvg{g_t} \right) \left( \Delta_\LB \varphi_i \right) \dd v, \quad i = 0.\dots,d, \\
    J_{ij}(t,x) := \int_{\S^{d-1}}
    \left[ \langle g_t \rangle - g_t \right] v_j  \varphi_i \dd v, \quad i=0,\dots, d, \ j=0,\dots, d,
  \end{dcases}
\end{align*}
which proves~\eqref{eq:claim-dec}-\eqref{eq:claim-am}. Going back to~\eqref{eq:inverse-div} we compute
\begin{align*}
  \int_0^{T^*} \int_{\Sigma} \left( \vAvg{g_t} - \gAvg{g} \right)^2 \dd x \dd t
  & = - \sum_{i=0} ^d \int_\sU \vec{F}_i K_i \dd x \dd t -  \sum_{i,j=0} ^d \int_\sU \vec{F}_i \partial_j J_{ij} \dd x \dd t \\
  & = - \sum_{i=0} ^d \int_\sU \vec{F}_i K_i \dd x \dd t + \sum_{i,j=0} ^d \int_\sU \partial_j \vec{F}_i J_{ij} \dd x \dd t
\end{align*}
where we have used the Dirichlet conditions again. Using the
$H^1(\sU)$ bound on $\vec{F}$ in~\eqref{eq:divergence}
and~\eqref{eq:claim-am} we deduce
\begin{align*}
  \int_0^{T^*} \int_{\Sigma} \left( \vAvg{g_t} - \gAvg{g} \right)^2 \dd x \dd t
  & \le \sqrt{d^2+d}\, \| \vec{F} \|_{H^1(\sU)} \left( \sum_{i=0} ^d \| K_i \|_{L^2(\sU)} ^2 + \sum_{i,j=0} ^d \|J_{ij}\|_{L^2(\sU)} ^2 \right)^{\frac12} \\
  & \le \sqrt{d^2+d}\, C_\sD \sqrt{C_3} \left\| \vAvg{g_t} - \gAvg{g} \right\|_{L^2(\sU)} \left( \int_0^{T^*} \cD(g_t) \dd t \right)^{\frac12}
\end{align*}
which implies by splitting the square
\begin{equation}
  \label{eq:local-avg-estim}
  \int_0^{T^*} \int_{\Sigma} \left( \vAvg{g_t} - \gAvg{g} \right)^2 \dd x \dd t \le C_\delta \int_0^{T^*} \cD(g_t) \dd t + \delta \int_0 ^{T^*} \| g_t \|^2 _{L^2(\T^d \times \S^{d-1})} \dd t
\end{equation}
for all $\delta >0$ and some corresponding constant $C_\delta$.

To finish the proof of~\eqref{eq:quant}, we need to estimate the
global average \(\gAvg{g}\) which we compare to the zero mass
condition up to error terms controlled by the dissipation. To relate
it to the zero mass condition~\eqref{eq:zeromass} introduce
\begin{equation}\label{eq:def-psi}
  \forall \, (t,z) \in [0,T^*] \times \T^d \times \S^{d-1}, \quad
  \psi(t,z) = \psi_t(z) := \frac{\chi(z)}{\int_0 ^{T^*} \chi(Z_{s-t}(z)) \dd s}.
\end{equation}
which is well-defined since the denominator is uniformly positive
thanks to~\eqref{eq:gcc}. The function $\psi$ is bounded in
$C^1([0,T^*] \times \T^d \times \S^{d-1})$, and satisfies
$\supp \psi = [0,T^*] \times \supp \sigma$, and, most importantly,
\begin{equation}
  \label{eq:avg-penal}
  \forall \, z \in \T^d \times \S^{d-1}, \quad \int_0 ^{T^*} \psi_t(Z_t(z)) \dd t = 1.
\end{equation}

By the conservation of mass, we find that \(h = 1 - \frac{1}{T^*} \vAvg{\psi}\)
has mass zero over \(\sU\). Hence we can apply \cref{lem:div} to
find \(\vec{F}\) with the properties of the lemma. We then find
\begin{equation*}
  \gAvg{g} - \frac{1}{T^*} \int_{[0,T^*] \times \T^d} \vAvg{g}(t,x)\,
  \vAvg{\psi}(t,x) \dd t \dd x
  = \int_{\sU} \vAvg{g}(t,x) \nabla \cdot \vec F \dd t \dd x.
\end{equation*}
Hence we can use~\eqref{eq:claim-dec}-\eqref{eq:claim-am} as before to find a constant \(C_4\) so that
\begin{equation}
  \label{eq:compare-averages}
  \left| \gAvg{g} -
    \frac{1}{T^*} \int_{[0,T^*] \times \T^d} \vAvg{g}(t,x)\,
      \vAvg{\psi}(t,x) \dd t \dd x
  \right|
  \le C_4 \left(\int_0^{T^*} \cD(g_t) \dd t\right)^{\frac 12}.
\end{equation}
We now estimate the $\psi$-weighted average as
\begin{align*}
  &\frac{1}{T^*} \int_{[0,T^*] \times \T^d} \vAvg{g}(t,x)\,
      \vAvg{\psi}(t,x) \dd t \dd x\\
  & = \int_{[0,T^*] \times \T^d \times \S^{d-1}} \vAvg{g}(t,x) \psi(t,x,v) \dd t \dd x \dd v \\
  & = \int_{[0,T^*] \times \T^d \times \S^{d-1}} \left[ M \vAvg{g_t}(x) - g_t(x,v) \right] \psi(t,x,v) \dd t \dd x \dd v \\
  & \qquad \qquad + \int_{[0,T^*] \times \T^d \times \S^{d-1}} g(t,x,v) \psi(t,x,v) \dd t \dd x \dd v =: J_1 + J_2.
\end{align*}

The first term $J_1$ is controlled by the
micro-coercivity~\eqref{eq:micro} and $\psi \lesssim \sigma$ for a
constant \(C_5\) as
\begin{equation}
  \label{eq:J1}
  J_1 \le C_5 \left( \int_{[0,T^*] \times \T^d \times \S^{d-1}} \sigma \left[ M \vAvg{g_t} - g_t \right]^2 \dd t \dd x \dd v \right)^{\frac12} \le C_5C_P \left( \int_0 ^{T^*} \cD(g_t) \dd t \right)^{\frac12}.
\end{equation}
We rewrite the second term $J_2$ by Duhamel's principle along the
transport flow as
\begin{align*}
  J_2 & = \int_{[0,T^*] \times \T^d \times \S^{d-1}} g_\init(Z_{-t}(z)) \psi_t(z) \dd t \dd z \\
      & \qquad \qquad + \int_{[0,T^*] \times
        \T^d \times \S^{d-1}} \int_0^t \sigma(X_{-(t-s)}(z))
        \Delta_\LB g_s(Z_{-(t-s)}(z)) \psi_t(z) \dd s \dd t \dd z \\
      & = \int_{\T^d \times \S^{d-1}} g_\init(z) \left( \int_0 ^{T^*} \psi_t(Z_t(z)) \dd t \right) \dd z \\
      & \qquad \qquad + \int_{[0,T^*] \times
        \T^d \times \S^{d-1}} \int_0^t \sigma(x) \Delta_\LB g_s(z)
        \psi_t(Z_{t-s}(z)) \dd s \dd t \dd z =: J_{21} + J_{22}
\end{align*}
and $J_{21}=0$ because of~\eqref{eq:avg-penal}
and~\eqref{eq:zeromass}, and the second term is estimated by
integration by parts:
\begin{align*}
  \left|J_{22}\right|
  & =
    \left| \int_{[0,T^*] \times \T^d \times \S^{d-1}} \int_0^t
    \sigma(x) \nabla_v g_s(z) \cdot \nabla_v \left[ \psi_t(x+(t-s)v,v)
    \right] \dd s \dd t \dd z \right| \\
  & \le C_6 \left( \int_0 ^{T^*} \cD(g_t) \dd t \right)^{\frac12}
\end{align*}
for some constant $C_6>0$. Together with~\eqref{eq:local-avg-estim}
and~\eqref{eq:J1} it concludes the proof of~\eqref{eq:quant}.

Let us finally extend the argument when \(\Sigma\) has finitely many connected components \(\Sigma_1,\dots,\Sigma_k\). For each $i = 1, \dots, k$, we define \(\sU_i = (0,T^*) \times \Sigma_i\) and
\begin{align*}
  \gAvg{g}_i := \frac{1}{m_i} \int_{[0,T^*] \times \Sigma_i} \vAvg{g}(t,x)
  \dd t \dd x
  \quad \mbox{ with } m_i & := |[0,T^*] \times \Sigma|.
\end{align*}

Arguing on each component as we did in the estimate~\eqref{eq:local-avg-estim}, we get
\begin{equation}
  \label{eq:local-avg-estim-comp}
  \int_0^{T^*} \int_{\Sigma_i} \left( \vAvg{g_t} - \gAvg{g}_i \right)^2 \dd x \dd t
  \le C_\delta \int_0^{T^*} \cD(g_t) \dd t + \delta \int_0 ^{T^*} \| g_t \|^2 _{L^2(\T^d \times \S^{d-1})} \dd t.
\end{equation}

We then prove, for each pair $i \not = j \in \{1,\dots,k \}$, that \(|\gAvg{g}_i - \gAvg{g}_j|^2\) is controlled by \(\int_0^{T^*} \cD(g_t) \dd t\). Indeed, all components are connected by the transport flow provided $T^*$ is chosen large enough (without loss of generality) so that there are smooth weights $w_i \ge 0$ over \((0,T^*) \times \T^d \times \S^{d-1}\) with unit masses and
\(\supp w_i \subset \sU_i \times \S^{d-1}\), and smooth compactly
supported \(\psi_{ij}=\psi_{ij}(t,x,v)\) solutions to
$\partial_t \psi_{ij} - v\cdot \nabla_x \psi_{ij} = w_i - w_j$.
Then integrating the equation on $g$ against $\psi_{ij}$ and using~\eqref{eq:micro} shows that
\begin{equation*}
  \left| \int_{[0,T^*] \times \T^d} \vAvg{g} \, \vAvg{w_i} \dd t \dd x
    - \int_{[0,T^*]   \times \T^d} \vAvg{g} \, \vAvg{w_j} \dd t \dd x\right|
  \lesssim \left(\int_0^{T^*} \cD(g_t) \dd t\right)^{\frac 12},
\end{equation*}
and arguing as in the proof of~\eqref{eq:compare-averages} we can
prove that each $\int_{[0,T^*] \times \T^d} \vAvg{g} \, \vAvg{w_i}\dd t \dd x$ is close to $\gAvg{g}_i$ up to an error of order $\left(\int_0^{T^*} \cD(g_t) \dd t\right)^{\frac 12}$.

We then construct \(\psi\) as in~\eqref{eq:def-psi} and, using
the zero mass condition~\eqref{eq:zeromass}, we can argue as above to
prove, for some constants \(\alpha_1,\dots,\alpha_k > 0\),
\begin{equation*}
  \alpha_1 \gAvg{g}_1 + \dots + \alpha_k \gAvg{g}_k
  \lesssim \left(\int_0^{T^*} \cD(g_t) \dd t\right)^{\frac 12}.
\end{equation*}
Together with the control on the differences between the $\gAvg{g}_i$'s this implies the result.

\section*{Acknowledgements}

All authors acknowledge partial support from the ERC grant MATKIT
grant. HD acknowledge the grant ANR-18-CE40-0027 of the French
National Research Agency (ANR).

\bibliographystyle{siam}
\bibliography{lit}

\end{document}